# Asymptotics of the solution of the Cauchy problem for a singularly perturbed system of hyperbolic equations[*]

**Nesterov A.V.** [1][0000-0002-4702-4777]

[1] PLEKHANOV Russian University of Economics, Stremyanny lane 36, Moscow, 117997, Russia
andrenesterov@yandex.ru

**Abstract.** An asymptotic expansion with respect to a small parameter of a singularly perturbed system of hyperbolic equations, describing vibrations of two rigidly connected strings is constructed. Under certain conditions imposed on these problems, the principal term of the asymptotic expansion of the solution is described by the Korteweg-de Vries equation.

## 1 Introduction

Consider a system of hyperbolic equations, describing the joint vibrations of two interconnected strings

$$\begin{cases} \varepsilon^3(u_{tt} - k_1 u_{xx}) = -au + bv + \varepsilon^2 f(u,v), \\ \varepsilon^3(v_{tt} - k_2 v_{xx}) = au - bv - \varepsilon^2 f(u,v), \end{cases} \quad (1)$$

with initial conditions

$$\begin{cases} u(x,0) = \overset{0}{u}(x), u_t(x,0) = \varphi(x), \\ v(x,0) = \frac{a}{b}\overset{0}{u}(x), v_t(x,0) = \frac{a}{b}\varphi(x), \end{cases} \quad (2)$$

where: $\{u(x,t), v(x,t)\}$ - solution, $\{x,t\} \in H = \{|x| < \infty; 0 \le t \le T, T > 0\}$; $0 < \varepsilon \ll 1$ small parameter, $k_1, k_2, a > 0, b > 0 - const$, $f(u,v) \subset C^\infty$, functions

---

[*] This research was performed in the framework of the state task in the field of scientific activity of the Ministry of Science and Higher Education of the Russian Federation, project "Development of the methodology and a software platform for the construction of digital twins, intellectual analysis and forecast of complex economic systems", grant no. FSSW-2020-0008.

$\overset{0}{u}(x), \varphi(x)$ have a certain smoothness (as discussed below) and satisfy the condition
$$\left|\overset{0}{u}(x), \varphi(x)\right| < Ce^{-\kappa x^2}, C, \kappa > 0.$$

Problem (1)-(2) is singularly perturbed in the critical case [1]. The study of the asymptotics of the solution of problem (1) - (2) is a continuation of papers [2]-[3], which investigated the asymptotics of solutions of singularly perturbed transport equations in the critical case. Problem (1)-(2) can be considered as one of the variants of generalization of the system of coupled oscillators, in which the Pasta-Ulam-Fermi paradox appears, which have been studied in many works, in particular, in [4],[5].

The asymptotic expansion (AE) with respect to a small parameter of the solution of problem (1)-(2) when the condition for initial data up to an arbitrary order $N$ is fulfilled, we will build in the form

$$\begin{pmatrix} u(x,t,\varepsilon) \\ v(x,t,\varepsilon) \end{pmatrix} = \begin{pmatrix} \bar{u}(x,t,\varepsilon) \\ \bar{v}(x,t,\varepsilon) \end{pmatrix} + \begin{pmatrix} Su(\zeta,\tau,\varepsilon) \\ Sv(\zeta,\tau,\varepsilon) \end{pmatrix} =$$
$$= \sum_{i=0}^{N} \varepsilon^i \left( \begin{pmatrix} \bar{u}_i(x,t) \\ \bar{v}_i(x,t) \end{pmatrix} + \begin{pmatrix} S_i u(\zeta,\tau) \\ S_i v(\zeta,\tau) \end{pmatrix} \right) + R_N = U_N + R_N \quad , \quad (3)$$

here $\bar{u}, \bar{v}$ is the regular part of AP, $Su, Sv$ are terms describing AE solutions along certain lines - "pseudo-characteristics" of system (1), $R$ is the residual term, $\tau = t/\varepsilon^2$, the type of variables $\zeta$ is described below.

In this paper, we restrict ourselves to constructing only the principal terms of AE (3).

In the case of smooth initial conditions $\overset{0}{u}(x), \varphi(x) \subset C^\infty_{[-\infty,+\infty]}$, the construction of AE (3) does not cause any special difficulties and is carried out by the method of boundary functions of A. B. Vasilyeva and V. F. Butuzov [1]. Of greater interest is the construction of AE in the case when the initial conditions are either piecewise continuous or have regions of large gradients.

## 2 Construction of a regular part of AE under smooth initial conditions.

For smooth initial conditions, the construction of the principal terms of the regular part of the AE is not difficult, but it is important for a more interesting case – initial conditions of the "burst" type.

The regular part of AE (3) has the form
$$\begin{pmatrix} \bar{u}(x,t,\varepsilon) \\ \bar{v}(x,t,\varepsilon) \end{pmatrix} = \sum_{i=0}^{3} \varepsilon^i \begin{pmatrix} \bar{u}_i(x,t) \\ \bar{v}_i(x,t) \end{pmatrix} + \dots \tag{4}$$

Terms with numbers $i = 1, 2, 3$ play an auxiliary role, while terms of a higher order of smallness are indicated by an ellipsis. Substitute AE (4) in the function $f(u,v)$ and represent it as

$$f(\bar{u}(x,t,\varepsilon),\bar{v}(x,t,\varepsilon)) =$$
$$= f(\bar{u}_0(x,t) + \varepsilon\bar{u}_1(x,t) + \varepsilon^2\bar{u}_2(x,t) + \varepsilon^3\bar{u}_3(x,t) + ...,$$
$$\bar{v}_0(x,t) + \varepsilon\bar{v}_1(x,t) + \varepsilon^2\bar{v}_2(x,t) + \varepsilon^3\bar{v}_3(x,t) + ...) = \quad (5)$$
$$= f(\bar{u}_0,\bar{v}_0) + \varepsilon F_1 + \varepsilon^2 F_2 + ...$$

where the terms $F_i$ depend on $\bar{u}_j, \bar{v}_j$ with numbers $j < i$.

Following [1], to determine the terms of expansion (4), we substitute expansion (4) and (5) in the system (1) and equate the left and right terms with the same powers of the parameter $\varepsilon$.

If the parameter degree is zero, we get

$$\varepsilon^0 : \begin{cases} -a\bar{u}_0 + b\bar{v}_0 = 0, \\ a\bar{u}_0 - b\bar{v}_0 = 0, \end{cases}$$

This system is solvable, hence

$$\bar{v}_0 = \frac{a}{b}\bar{u}_0 \quad (6)$$

where $\bar{u}_0(x,t)$ - is an yet undefined function.

For the first degree of the parameter, we get

$$\varepsilon^1 : \begin{cases} -a\bar{u}_1 + b\bar{v}_1 = f(\bar{u}_0,\bar{v}_0), \\ a\bar{u}_1 - b\bar{v}_1 = -f(\bar{u}_0,\bar{v}_0), \end{cases}$$

This system is solvable, hence

$$\bar{v}_1 = \frac{a}{b}\bar{u}_1 + \frac{1}{b}f(\bar{u}_0,\bar{v}_0),$$

where $\bar{u}_1(x,t)$ - is an yet undefined function.

For the second power of the parameter, we get

$$\varepsilon^2 : \begin{cases} -a\bar{u}_2 + b\bar{v}_2 = F_1, \\ a\bar{u}_1 - b\bar{v}_1 = -F_1, \end{cases}$$

This system is solvable, hence

$$\bar{v}_2 = \frac{a}{b}\bar{u}_2 + \frac{1}{b}F_1,$$

where $\bar{u}_1(x,t)$ - is an yet undefined function.

For the third power of a small parameter, we obtain

$$\varepsilon^3 : \begin{cases} -a\bar{u}_3 + b\bar{v}_3 = \bar{u}_{0,tt} - k_1\bar{u}_{0,xx} + F_2, \\ a\bar{u}_3 - b\bar{v}_3 = \bar{v}_{0,tt} - k_2\bar{v}_{0,xx} - F_2, \end{cases} \quad (7)$$

For the solvability of the system of equations (7), it is necessary (and sufficient) to satisfy the condition

$$\bar{u}_{0,tt} - k_1\bar{u}_{0,xx} + \bar{v}_{0,tt} - k_2\bar{v}_{0,xx} = 0 \tag{8}$$

By substituting (6) in (8), we obtain the equation for determining $\bar{u}_0$.

Let's introduce the notation

$$\frac{bk_1 + ak_2}{a+b} = k, \min(k_1,k_2) < k < \max(k_1,k_2) \tag{9}$$

In this notation, the equation for the definition $\bar{u}_0$ takes the form

$$\bar{u}_{0,tt} - k\bar{u}_{0,xx} = 0. \tag{10}$$

The function $\bar{v}_0$ is determined from the algebraic relation

$$\bar{v}_0 = \frac{a}{b}\bar{u}_0. \tag{11}$$

Let's set the initial conditions for determining $\bar{u}_0$

$$\bar{u}_0(x,0) = \overset{0}{u}(x), \bar{u}_{0,t}(x,0) = \varphi(x). \tag{12}$$

Due to the special choice of initial conditions and relation (11), the initial conditions for $v_0$ $v_0(x,0) = \frac{a}{b}\overset{0}{u}(x), v_{0,t}(x,0) = \frac{a}{b}\varphi(x)$ are automatically satisfied. Boundary functions describing fast transients in a small neighborhood of the line $t = 0$, are absent.

This completes the construction of the principal terms of AP (4).

*Theorem* (on the evaluation of the residual term of AE (4)).

If the conditions $f(u,v) \subset C^\infty$, $\overset{0}{u}(x), \varphi(x) \subset C^\infty_{(-\infty,+\infty)}$ $\left|\overset{0}{u}(x), \varphi(x)\right| < Ce^{-\kappa x^2}, C, \kappa > 0$, are met, the solution of problem (1)-(2) can be represented in the form $\begin{pmatrix} \bar{u}(x,t,\varepsilon) \\ \bar{v}(x,t,\varepsilon) \end{pmatrix} = \begin{pmatrix} \bar{u}_0(x,t) \\ \bar{v}_0(x,t) \end{pmatrix} + \begin{pmatrix} Ru \\ Rv \end{pmatrix} = U + R$, $(\bar{u}_0, \bar{v}_0)$ defined above, and $R$ there is a solution to the problem

$$\begin{cases} \varepsilon^3(Ru_{tt} - k_1 Ru_{xx}) = -aRu + bRv + \varepsilon^2 r(Ru,Rv) + r_1, \\ \varepsilon^3(Rv_{tt} - k_2 Rv_{xx}) = aRu - bRv - \varepsilon^2 r(Ru,Rv) + r_2, \\ Ru(x,0) = 0, Ru_t(x,0) = 0, \\ Rv(x,0) = 0, Rv_t(x,0) = 0. \end{cases}$$

Here $r(Ru,Rv) = f(\bar{u},\bar{v}) + (f(\bar{u} + Ru, \bar{v} + Rv), r_1, r_2 = O(\varepsilon^2)$. The proof of the theorem follows from the function construction algorithm $(\bar{u}_0, \bar{v}_0)$.

# 3 Building an AE under initial conditions of the "splash" type

More interesting is the case when the initial conditions contain regions of asymptotically fast changes, for example, of the "smoothed step" or "surge" type.

Let the initial conditions for equation (1) have the form

$$\begin{cases} u(x,0) = \overset{0}{u}(\frac{x}{\varepsilon}), u_t(x,0) = \varphi(\frac{x}{\varepsilon}), \\ v(x,0) = \frac{a}{b}\overset{0}{u}(\frac{x}{\varepsilon}), v_t(x,0) = \frac{a}{b}\varphi(\frac{x}{\varepsilon}), \end{cases} \quad (13)$$

Here functions $\overset{0}{u}(x), \varphi(x)$ have a certain smoothness and satisfy the condition $\left| \overset{0}{u}(x), \varphi(x) \right| < Ce^{-\kappa x^2}, C, \kappa > 0.$ Condition (13) specify the initial data in the form of an asymptotically narrow "cap", "splash".

We will build for AE solutions of system (1) with initial conditions (13) in the form

$$\begin{pmatrix} u(x,t,\varepsilon) \\ v(x,t,\varepsilon) \end{pmatrix} = \begin{pmatrix} \bar{u}(x,t,\varepsilon) \\ \bar{v}(x,t,\varepsilon) \end{pmatrix} + \begin{pmatrix} Su(\zeta,\tau,\varepsilon) \\ Sv(\zeta,\tau,\varepsilon) \end{pmatrix} + R =$$

$$= \sum_{i=0}^{3} \varepsilon^i \left( \begin{pmatrix} \bar{u}_i(x,t) \\ \bar{v}_i(x,t) \end{pmatrix} + \begin{pmatrix} S_i u(\zeta,\tau) \\ S_i v(\zeta,\tau) \end{pmatrix} \right) + R_3 = U_3 + R_3$$

The regular part of the AP solution $\bar{U}_N$ under such initial conditions is obviously zero and the AE takes the form

$$\begin{pmatrix} u(x,t,\varepsilon) \\ v(x,t,\varepsilon) \end{pmatrix} = \sum_{i=0}^{3} \varepsilon^i \begin{pmatrix} S_i u(\zeta,\tau) \\ S_i v(\zeta,\tau) \end{pmatrix} + R_3 \quad (14)$$

Accordingly, the functions $Su, Sv$ must satisfy system (1) and the initial conditions (13) or part of these conditions. Defining new stretched variables

$$\zeta_{1,2} = \varepsilon^{-1}(x \pm kt), \quad (15)$$

where $k$ is defined by expression (9). Given (14), we transform $f(Su, Sv)$

$$f(Su, Sv) = f(\sum_{i=0}^{3} \varepsilon^i Su_i + R_3 u, \sum_{i=0}^{3} \varepsilon^i Sv_i + R_3 u) =$$
$$= f(Su_0, Sv_0) + \varepsilon F_1 + \varepsilon^2 F_2 + ... \quad (16)$$

where the terms $F_i$ depend on $S_j u, S_j v$ with numbers $j < i$.

Building the function

$$S^I(\zeta_1, t, \varepsilon) = \begin{pmatrix} S^I u(\zeta_1, t, \varepsilon) \\ S^I v(\zeta_1, t, \varepsilon) \end{pmatrix} = \sum_{i=0}^{2} \varepsilon^i \begin{pmatrix} S^I_i u(\zeta_1, t) \\ S^I_i v(\zeta_1, t) \end{pmatrix} + R, \quad (17)$$

where $\zeta_1 = \dfrac{x - kt}{\varepsilon}$ ( the index $I$ of the functions $S$ and the index $1$ of the variable $\zeta$ omitted are below). Functions with numbers 1,2 in decomposition (17) play an auxiliary role.

Moving in system (1) from variables $(x,t)$ to variables $(\zeta,t)$

$$\begin{cases} \varepsilon^3 u_{tt} - 2k\varepsilon^2 u_{\zeta t} + \varepsilon(k^2 - k^2_{\ 1})u_{\zeta\zeta} = -au + bv + \varepsilon f(u,v), \\ \varepsilon^3 v_{tt} - 2k\varepsilon^2 v_{\zeta t} + (k^2 - k^2_{\ 2})\varepsilon v_{\zeta\zeta} = au - bv - \varepsilon f(u,v), \end{cases}$$

let's rewrite it as

$$\begin{cases} -au + bv = (\varepsilon^3 u_{tt} - 2k\varepsilon^2 u_{\zeta t} + \varepsilon(k^2 - k^2_{\ 1})u_{\zeta\zeta}) - \varepsilon f(u,v), \\ au - bv = (\varepsilon^3 v_{tt} - 2k\varepsilon^2 v_{\zeta t} + (k^2 - k^2_{\ 2})\varepsilon v_{\zeta\zeta}) + \varepsilon f(u,v), \end{cases} \quad (18)$$

Substituting expansions (14) and (17) into (18), equating the coefficients for the same degrees $\varepsilon$ in the left and right sides of the equations, we obtain equations for determining the terms of the expansion (14).

If the parameter degree is zero, we get

$$\varepsilon^0 : \begin{cases} -aS_0 u + bS_0 v = 0, \\ aS_0 u - bS_0 v = 0, \end{cases}$$

This system is solvable, hence

$$S_0 v = \dfrac{a}{b} S_0 u, \quad (19)$$

where $S_0 u(\zeta, t)$ - is an yet undefined function.

For the first degree of the parameter, we get

$$\varepsilon^1 : \begin{cases} -aS_1 u + bS_1 v = (k^2 - k_1^{\ 2})S_0 u_{\zeta\zeta} - f(S_0 u, S_0 v), \\ aS_1 u - bS_1 v = (k^2 - k_2^{\ 2})S_0 v_{\zeta\zeta} + f(S_0 u, S_0 v). \end{cases} \quad (20)$$

This system is solvable because

$$((k^2 - k_1^{\ 2})S_0 u_{\zeta\zeta} - f(S_0 u, S_0 v)) + (k^2 - k_2^{\ 2})S_0 v_{\zeta\zeta} + f(S_0 u, S_0 v) =$$
$$= ((k^2 - k_1^{\ 2}) + \dfrac{a}{b}(k^2 - k_2^{\ 2}))S_0 u_{\zeta\zeta} = \dfrac{a(k^2 - k_2^{\ 2}) + b(k^2 - k_1^{\ 2})}{b} S_0 u_{\zeta\zeta} = 0$$

by virtue of definition $k^2$ (9).
From (19) and (20) we obtain

$$S_1 v = \frac{a}{b} S_1 u + \frac{1}{b}((k^2 - k_1^2) S_0 u_{\zeta\zeta} - f(S_0 u, \frac{a}{b} S_0 u)), \qquad (21)$$

where $S_1 u(\zeta, t)$ - is an yet undefined function.

For the second power of the parameter, we get

$$\varepsilon^2 : \begin{cases} -a S_2 u + b S_2 v = (k^2 - k_1^2) S_1 u_{\zeta\zeta} - 2k S_0 u_{\zeta t} - F_1, \\ a S_2 u - S_2 v = (k^2 - k_2^2) S_1 v_{\zeta\zeta} - 2k S_0 v_{\zeta t} + F_1. \end{cases}$$

The solvability condition of this system has the form

$$(k^2 - k_1^2) S_1 u_{\zeta\zeta} - 2k S_0 u_{\zeta t} + (k^2 - k_2^2) S_1 v_{\zeta\zeta} - 2k S_0 v_{\zeta t} = 0 \qquad (22)$$

Excluding from (22) $S_1 v(\zeta, t)$ by (21), we obtain the equation for the definition $S_0 u(\zeta, t)$.

Let's introduce the notation

$$\frac{(k^2 - k_2^2)(k^2 - k_1^2)}{2k(a+b)} = K,$$

$$-\frac{b}{2k(a+b)}(k^2 - k_2^2) f(S_0 u, \frac{a}{b} S_0 u) = h(S_0 u). \qquad (23)$$

Taking into account the notation (23), the equation for determining $S_0 u(\zeta, t)$ is written as

$$(-S_0 u_t + K S_0 u_{\zeta\zeta\zeta} - h(S_0 u)_\zeta)_\zeta = 0$$

or

$$-S_0 u_t + K S_0 u_{\zeta\zeta\zeta} - h(S_0 u)_\zeta = \Psi(t),$$

where $\Psi(t)$ is an arbitrary function. Based on the requirement $S_0 u \underset{\zeta \to \pm\infty}{\to} 0$, it is natural to put $\Psi(t) = 0$. This gives the final form of the equation for the function $S_0 u$

$$-S_0 u_t + K S_0 u_{\zeta\zeta\zeta} - h(S_0 u)_\zeta = 0. \qquad (24)$$

The functions $S^{II}(\zeta_2, t, \varepsilon)$

$$S^{II}(\zeta_2, \tau, \varepsilon) = \begin{pmatrix} S^{II} u(\zeta_2, \tau, \varepsilon) \\ S^{II} v(\zeta_2, \tau, \varepsilon) \end{pmatrix} = \sum_{i=0}^{3} \varepsilon^i \begin{pmatrix} S^{II}_i u(\zeta_2, \tau) \\ S^{II}_i v(\zeta_2, \tau) \end{pmatrix} + R \qquad (25)$$

are constructed similarly.

where $\zeta_2 = \dfrac{x+kt}{\varepsilon}$ ( indexes $II$ $S^{II}$ and $2$ variable $\zeta_2$ can be omitted below).

Moving from variables $(x,t)$ to variables $(\zeta_2,t)$ in system (1) (the index $2$ can be omitted below), we reduce it to the form

$$\begin{cases} \varepsilon^3(u_{tt} + 2k\varepsilon^{-1}u_{\zeta t} + k^2\varepsilon^{-2}u_{\zeta\zeta} - k^2_1\varepsilon^{-2}u_{\zeta\zeta}) = -au + bv + \varepsilon^2 f(u,v), \\ \varepsilon^3(v_{tt} + 2k\varepsilon^{-1}v_{\zeta t} + k^2\varepsilon^{-2}v_{\zeta\zeta} - k^2_2\varepsilon^{-2}v_{\zeta\zeta}) = au - bv - \varepsilon^2 f(u,v), \end{cases} \quad (26)$$

Substituting (25) and (16) in the system (26), by analogy, we obtain the equation for determining $S^{II}_0(\zeta_2,t)$

$$S^{II}_0 u_t + K S^{II}_0 u_{\zeta\zeta\zeta} - h(S^{II}_0 u)_\zeta = 0, \qquad (27)$$

and the ratio to determine $S^{II}_0 v$ : $S^{II}_0 v = \dfrac{a}{b} S^{II}_0 u$.

As a result, the main terms of the AE solution of problem (1), (13) are constructed in the form

$$\begin{pmatrix} u(x,t,\varepsilon) \\ v(x,t,\varepsilon) \end{pmatrix} = \begin{pmatrix} S^I_0 u(\zeta_1,t) + S^{II}_0 u(\zeta_2,t) \\ S^I_0 v(\zeta_1,t) + S^{II}_0 v(\zeta_2,t) \end{pmatrix} + R \qquad (28)$$

These variables $(\zeta_1,\zeta_2)$ are defined by relations (15), and the functions $S^I_0 u(\zeta_1,t), S^{II}_0 u(\zeta_2,t), S^I_0 v(\zeta_1,t), S^{II}_0 v(\zeta_2,t)$ satisfy equations (24), (27).

To fully define the functions $S^I_0 u(\zeta_1,t), S^{II}_0 u(\zeta_2,t), S^I_0 v(\zeta_1,t), S^{II}_0 v(\zeta_2,t)$, you need to set the initial conditions, which will be the subject of further research.

## Conclusion.

The system of equations (1) can be considered as one of the variants of generalization of the problem of coupled oscillators, considered in many works, in particular, in [4],[5], in which the phenomenon of the Pasta-Ulama-Fermi paradox is analyzed and the connection of the system of coupled oscillators with the Korteweg-de Vries equation is noted. In this formulation, the AE solutions of system (1) with initial conditions of the burst type (13) in the first approximation have the form (28). $S^I_0 u, S^{II}_0 u$ are described by equations (24), (27) – generalized Korteweg -de Vries equations. In the continuation of this study, we plan to determine the initial conditions for equations (23), (27). Additional boundary functions will need to be introduced in order to satisfy the initial conditions (13), even in the first approximation in AE (28).